\documentclass{amsart}

\usepackage{amssymb}
\usepackage[all]{xy}
\usepackage{hyperref}

\usepackage{enumitem}   

\setlist[enumerate]{itemsep=.2em,topsep=.2em,leftmargin=1.25em,itemindent=2.0em}

\newtheorem{thm}{Theorem}

\newtheorem{lem}[thm]{Lemma}

\theoremstyle{definition}

\newtheorem{say}[thm]{}
\newtheorem{exmp}[thm]{Example}

\newtheorem{rem}[thm]{Remark}

\newtheorem*{ack}{Acknowledgments}      

\newtheorem{defn-thm}[thm]{Definition--Theorem}  
\newtheorem{defn-lem}[thm]{Definition--Lemma}  

\theoremstyle{remark}

\renewcommand{\c}[0]{{\mathbb C}}  

\renewcommand{\o}[0]{{\mathcal O}} 

\renewcommand{\a}[0]{{\mathbb A}}

\newcommand{\p}[0]{{\mathbb P}}
\newcommand{\f}[0]{{\mathbb F}}

\newcommand{\mult}[0]{\operatorname{mult}}

\newcommand{\aut}[0]{\operatorname{Aut}}

\def\loccoh#1.#2.#3.#4.{H^{#1}_{#2}(#3,#4)}

\DeclareMathAlphabet{\mathchanc}{OT1}{pzc}%
                                {m}{it}

\usepackage[all]{xy}\xyoption{dvips}

\begin{document}
\bibliographystyle{amsalpha}


\title[Planar, rational curves]
      {Planar, rational curves   over $\f_2$, \\  whose only singularity is a double point}
 \author{J\'anos Koll\'ar}

 \begin{abstract}
   We exhibit planar, rational curves of large degree
   over $\f_2$ that have a unique singular point, which has multiplicity 2. In characteristic 0 such curves exist only for degrees up to $6$.
 \end{abstract}

 \maketitle

 The aim  of this note is to construct rational curves
 $C_d\subset \f_2\p^2$ of degree $d=2^n+2$ or $d=2^n$ that have a unique singular point, which has multiplicity 2; see Example~\ref{main.2.exmp}.

 In characteristic 0 such curves exist only for degrees up to $6$,
 see    Lemma~\ref{gzn.lem} and Remark~\ref{yosh.rem}.

 For $d\geq 8$ the curves $C_d$ cannot be lifted to characteristic 0
 in a way that preserves the standard resolution-of-singularities  blow-up sequence; see Remark~\ref{budef.rem}.
 This shows that the program outlined in
 \cite{ishii2025liftingsidealspositivecharacteristic} needs changes in higher dimensions.
 However, these curves do not seem to   
 shed light on whether the set of 
 log canonical thresholds depends on the  characteristic, which is one of  the main topics  of \cite{ishii2025liftingsidealspositivecharacteristic, ishii2026liftingsidealspositivecharacteristic-v4};
see Remark~\ref{lct.rem}.

 A general, planar,  rational curve $B_d\subset \p^2$ of degree $d$ has
 $\binom{d-1}{2}$ nodes, and it is not easy to obtain examples with few singular points. There seem to be very few curves with only one singular point, which is  a cusp.
These appear  prominently in the 
 Abhyankar-Moh-Suzuki theorem describing $\aut(\c^2)$
\cite{MR0338423, MR0379502}. 
  Those with a single  Puiseux pair are classified in  \cite{MR2280130}. 

It is especially restrictive to have only  cusps with multiplicity 2
   \cite{MR1016092, MR1802321}; these can give surprising examples of real polynomials \cite{MR2795944}.

   In positive characteristic, Artin-Schreier--type polynomials have been used to give non-standard embeddings of $\a^1$ into $\a^2$.
   This was most likely  observed by Abhyankar;  
\cite{MR409472} refers to them as  `known counterexamples.' 
   The curves  $C_d$ can be  given by Artin-Schreier--type equations as in (\ref{main.exmp.lem}.5),
 but their properties are  easier to see using an   explicit  parametrization.

\begin{exmp}\label{main.exmp}
  Fix a characteristic $p>0$. Let $q$ be a $p$-power, and $g(u)\in \f_p[u]$ a polynomial.
  Let $C^\circ_{q,g}\subset \a^2$ be the the image of the map
   $$
   \phi_{q,g}: t\mapsto  \bigl(t^q, g(t^p)+t\bigr)=:(x,y),
   $$
   and $C_{q,g}\subset \p^2$ its closure.
   \end{exmp}

 \begin{lem} \label{main.exmp.lem}
  For the curves $C_{q,g}$ in Example~\ref{main.exmp} the following hold.
   \begin{enumerate}
   \item  $\phi_{q,g}: \a^1\to \a^2$ is a closed embedding,
\item $C^\circ_{q,g}$ is a smooth, rational curve, 
     and
     \item $C_{q,g}$ has a  single point  at infinity. If $r:=p\cdot\deg g -q$ is positive, it has multiplicity $r$.
   \end{enumerate}
   \end{lem}

 Proof. The derivatve of $\phi_{q,g}$ is $(0,1)$, so
 it is an immersion.  Since $t^q=x$, it is also purely inseparable,
 hence an isomorphism.
 The equation of $C^\circ_{q,g}$ is obtained using
 $$
 y^q=g(t^{qp})+t^q=g(x^p)+x.
 \eqno{(\ref{main.exmp.lem}.4)}
 $$
 If $r:=p\cdot\deg g -q$ is positive, then
 the unique point at infinity is $(0{:}1{:}0)$, and the equation there is
   $$
 z^r=g(x^p)+xz^{q+r-1}.
 \qed
  \eqno{(\ref{main.exmp.lem}.5)}
     $$

   \begin{exmp} \label{main.2.exmp} Setting $p=2, q=2^n$  and $g(u)=u^{2^{n-1}+1}$, we get
     $$
     C_{2^n,2}:=\bigl(z^2y^{2^n}=x^{2^n+2}+xz^{2^n+1}\bigr)\subset \p^2.
     $$
     It is a rational curve of degree $2^n+2$.
     It has a unique singular point, which is a cusp of  multiplicity 2.

     It is thus of type  $A_m$, where $m=2^n(2^n+1)$.
     The latter is roughly  $(\deg C_{2^n,2})^2$.

     Note that for $m$ even, there are non-equivalent singularities of type $A_m$ in characteristic 2  with the same embedded resolution graphs. In the notation of 
     \cite{MR818299, MR1033443}, the normal form of the singularity of $C_{2^n,2}$ is $z^2+zx^r+x^{m+1}$ where $m=2^n(2^n+1)$ and
     $r=\frac12m+2^{n-1}+1$.

     Setting $p=2, q=2^n$  and $g(u)=u^{2^{n-1}-1}$, we get
     $$
     C_{2^n}:=\bigl(y^{2^n}=z^{2}+z^{2^n-1}\bigr)\subset \p^2.
     $$
     It is a rational curve of degree $2^n$.
     It has a unique singular point of type  $A_m$, where $m=(2^n-1)(2^n-2)$.
   \end{exmp}

   \begin{lem}\cite{MR1802321} \label{gzn.lem}  Let $B\subset \c\p^2$ be a reduced,  (not necessarily rational) curve of degree $2d$ with a  singularity of type $A_m$.   Then
     $ m\leq 3d(d-1)+1$,
     which is roughly  $\frac34(\deg B)^2$.
   \end{lem}

   Hint of proof. Let $S\to \p^2$ be the double cover  ramified along $B$.
   It has a surface singularity of type $A_m$. Thus the cohomology of the vanishing cycles of the $A_m$ surface singularity embeds into    $H^{2,2}(S_t, \c)$ of a
   smoothing $S_t$ of $S$. \qed
   \medskip

   \begin{rem}
   The upper bound in Lemma~\ref{gzn.lem} is not expected to be sharp for $d>6$.  The reducible curves
   $$
  D_{2d}:= \bigl((y-x^d)^2=y^{2d}\bigr)
  $$
  have a singularity of type  $A_{2d^2-1}$.
  This corresponds to 
  $m=\frac{1}{2}(\deg D_{2d})^2-1$.

  It is not easy to find curves $B$ with an $A_m$ singularity such that 
  $m>\frac{1}{2}(\deg B)^2$. 
       The largest known examples give asymptotically
       $m\approx \frac{7}{12}(\deg B)^2$ \cite[Sec.5]{orevkov-2012}.
\end{rem}

   \begin{rem}\label{yosh.rem}
     For $d\leq 6$ 
     \cite{MR533711} gives degree $d$, rational,  plane curves   $Y_d\subset \c\p^2$ with a single singular point, which has multiplicity 2.
     The curve $Y_6$ has an 
     $A_{19}$ singularity; the largest possible.
     Simpler  equations are given in \cite[5.8]{MR1900779}.  
   
    The curve   $ C_{4,2}$ has degree $6$ and an   $A_{20}$ singularity.
     It is possible that  $Y_6$ specializes to  $C_{4,2}$.
     Note that the singularity type changes from $A_{19}$ to $A_{20}$, but they both need 10 blow-ups to resolve.
     So the embedded resolution of $C_{4,2}$ as in (\ref{budef.rem}.1) might  be liftable to an embedded resolution of $Y_6$.
 \end{rem}

     \begin{rem}[Connection with  \cite{ishii2025liftingsidealspositivecharacteristic}]\label{budef.rem}
       Let $0\in X\subset \a^n$ be a singularity. Assume for simplicity
       that it can be resolved by a sequence of blowing up  points
       $$
       (Y_r, X_r, p_r) \to   \cdots\to(Y_1, X_1, p_1) \to  (Y_0, X_0, p_0)=(\a^n,X, 0),
       \eqno{(\ref{budef.rem}.1)}
       $$
       where $X_i\subset Y_i$  is   the birational transform of $X$, and
       $p_{i+1}$ lies over $p_i$.  Let $\pi_i:Y_i\to Y_0$ denote the composite.

       The discrepancies of the exceptional divisors with respect to the ideal sheaf of $X$ are determined by
       the numbers $\mult_{p_i}\pi_i^{-1}(p_0)$ and
       $\mult_{p_i} X_i$.

       Thus if we deform the blow-up sequence (\ref{budef.rem}.1)
       such that these multiplicities stay constant, then
       the discrepancies of the $\pi_r$-exceptional divisors are invariant under this deformation.
       Call these {\it discrepancy preserving} deformations of (\ref{budef.rem}.1).

       \cite{ishii2025liftingsidealspositivecharacteristic} notes that every exceptional divisor over $\a^n$ can be obtained as a
       $\pi_r$-exceptional divisor for a (somewhat more general)
       blow-up sequence as in  (\ref{budef.rem}.1).
       
       This is especially useful if $0\in X\subset \a^n$ is in characteristic $p$, and the general fiber of such a deformation is in characteristic $0$. It seems that in many cases there are discrepancy preserving liftings from characteristic $p$ to  characteristic $0$.
       The 2-dimensional cases are treated in \cite{MR584440, ishii2026liftingsidealspositivecharacteristic-v4}.

However, the curves in Example~\ref{main.2.exmp}  show that this  way of preserving 
       discrepancies in liftings is not always possible in higher dimensions.

       To get such examples, 
       let  $C_d=\bigl(g(x,y,z)=0\bigr)$ be a curve of degree $d$ in characteristic $2$ that has a single $A_{2r}$ singularity. For odd $m\gg 1$ consider
       $$
       X_g:=\bigl(g(x,y,z)+x^m+y^m+z^m=0\bigr)\subset \a^3.
       $$
       This can be resolved by first blowing up the origin, and then   $r$-times
       the singular point of the curve $C_d$.

       The resulting blow-up sequence has a discrepancy preserving lifting to
       characteristic $0$ iff the curve $C_d$ lifts to a curve $C'_d$ in
       characteristic $0$ that  has an $A_{2r}$ or $A_{2r-1}$ singularity.

       The smallest counterexample is
  $C_{8}$, which has degree $8$ and an   $A_{42}$ singularity.
    A reduced, octic curve $B_{8}\subset \c\p^2$ has at most an
     $A_{37}$ singularity by Lemma~\ref{gzn.lem}.
       
     The next one is     $C_{8,2}$, which has degree $10$ and an   $A_{72}$ singularity.
     A reduced, decic curve $B_{10}\subset \c\p^2$ has at most an
     $A_{61}$ singularity by Lemma~\ref{gzn.lem}.
     \end{rem}

     \begin{rem}[Log canonical thresholds]\label{lct.rem}
       The program of Ishii  aims to understand 
 log canonical thresholds  in positive  characteristic
 \cite{ishii2025liftingsidealspositivecharacteristic}.
 
     The log canonical threshold of the curve
     $C_{2^n, 2}\subset \f_2\p^2$ is $\frac12+\frac1{m+1}$, where $m=2^{2n}+2^n$.
     For $n\geq 2$ there is no curve $B\subset \c\p^2$ with the same degree and same log canonical threshold as $C_{2^n, 2}$.
     
     However,   the 
     log canonical threshold of the surface $X(g)\subset \a^3$
constructed from $C_{2^n, 2}$
is $3/(2^n+2)$. It is computed by the  exceptional divisor of the first blow-up, and  a general singular point with the same multiplicity has the same log canonical threshold.
 Thus these curves do not seem to  give further information about  the program in 
 \cite{ishii2025liftingsidealspositivecharacteristic}.
      \end{rem}

     Note,
        however, that the above non-liftable examples seem rather exceptional, and possibly their discrepancies can be handled in different ways. The next examples show that liftability may frequently hold.

     \begin{exmp}[Quartic surfaces] \label{quart.rem}
       I first tried   to look for non-liftable examples among K3 surfaces.
In  characteristic $0$ the rank of the N\'eron-Severi group of a smooth K3 is at most 20,
but in  characteristic $p$ it is  22 for supersingular ones.

So if $S$ is a supersingular K3 surface with Picard number 1
in characteristic $p$, then 
 there is no  lifting 
to  characteristic $0$ preserving {\em all} the singularity types.

Most  supersingular K3 surfaces with Picard number 1 have many singularities, and individually each singularity usually lifts, but there are a few examples where lifting even {\em one} singularity is not possible.
I thank Shimada and Sch\"utt for explaining these to me.

Working with the lists of \cite{MR2059747},
we discuss 4 extremal cases, but none are good enough to be used in Remark~\ref{budef.rem}.

\medskip

(\ref{quart.rem}.1) There is a  degree 4  K3  surface $S_1$ in  characteristic $2$ with a $D_{21}$ singularity. However, 
$S_1$  is not  a quartic in $\p^3$, but  a double cover of the quadric cone. Thus $S_1$ is naturally an octic surface in  $\p(1,1,2,4)$. So there is  an octic  surface in $\p(1,1,2,4)$
with a $D_{21}$ singularity in  characteristic $2$, but there is no such in char 0.

It seems likely that if $(s,S)$ is any 
quartic surface  in characteristic $p$ with with an ADE singularity at $s$,
then there is a lifting 
to  characteristic $0$ preserving the singularity type at $s$.

(\ref{quart.rem}.2) There is a K3  surface $S_2$ in  characteristic $11$ with  an $A_{21}$ singularity, but it has degree $22$.

(\ref{quart.rem}.3) There is a quasi-elliptic K3  surface $S_3$ in  characteristic $2$ with a $D_{20}$ singularity, and a section along which the surface is smooth.
 The Weierstrass equation realizes it is a surface in
 the projectivization of $\o_{\p^1}+\o_{\p^1}(8)+\o_{\p^1}(12)$, but there is no similar surface in  characteristic $0$.

 (\ref{quart.rem}.4) There is an elliptic K3  surface $S_4$ in  characteristic $2$ with  $A_{17}+3A_1$ singularities, and a section along which the surface is smooth.
 As before, there is no lifting that keeps all singularities, but
there is an elliptic K3  surface  in  characteristic $0$ with  $A_{17}+A_1$ singularities. 
\end{exmp}

The curves $C_{2^n,2}$ also lead to supersingular double planes.

\begin{say}[Supersingular double planes]\label{ss.dp.say}
  Over $\f_2$ consider the surfaces
$$
  S_{r}:=\bigl(u^2+u(zy^{2r}+x^{2r+1})+xz^{4r+1}+z^2y^{4r}+x^{4r+2}=0\bigr)\subset \p(1,1,1,2r+1).
  \eqno{(\ref{ss.dp.say}.1)}
  $$
 In the $y=1$ chart, the equation is
  $$
 g:=u^2+u(z+x^{2r+1})+xz^{4r+1}+z^2+x^{4r+2}=0.
 \eqno{(\ref{ss.dp.say}.2)}
    $$
    Thus
    $
    g_u=z+x^{2r+1},\
    g_z=u+xz^{4r}, \
    g_x=ux^{2r}+z^{4r+1}.
    $
   The first 2 define the monomial curve $M_r$ which is the image of
    $$
    \phi_r: t\mapsto \bigl(t, t^{2r+1}, t^{8r^2+4r+1}\bigr),
    $$
    and $g_x$ vanishes identically along $M_r$.
    Thus
    $$
    \begin{array}{rcl}
    \f_2[x,z,u]/(g, g_x, g_z, g_u)&\cong&
    \f_2[t]/\bigl(g(t, t^{2r+1}, t^{8r^2+4r+1})\bigr)\\
    &\cong&
    \f_2[t]/\bigl(t^{8r^2+6r+2}(t^{8r^2+2r}+1)\bigr).
    \end{array}
    $$
    Thus $S_r$ has an $A_{8r^2+6r+1}$ singularity at $t=0$.
For the other singularities, write  $2r=qr'$ where $q$ is a $2$-power and $r'$ is odd. Then we have 
$A_{q-1}$ singularities when $t$ is a $(4rr'+r')^{th}$ root of unity.
There are no other singular points.

If $r$ is a $2$-power, 
then the minimal resolution  $\bar S_r\to S_r$ contains 
$ 8r^2+6r+1+(4r+1)(2r-1)=16r^2+4r$ exceptional curves.

    Over the line $(z=0)$ the equation  (\ref{ss.dp.say}.2) becomes
    $u^2+ux^{2r+1}+x^{4r+2}=0$. Setting   $v:=ux^{-2r-1}$ we get
    $v^2+v+1=0$. If $\alpha_1, \alpha_2$ are its roots,
    then the preimage of $(z=0)$ has 2 components
    $$
    B_i:=\bigl(u+\alpha_ix^{2r+1}=z=0\bigr).
    $$
    Thus the Picard number of $\bar S_r$ is at least
    $16r^2+4r+2$.

    This equals the second Betti bumber of a double plane ramified along a smooth curve of degree $4r+2$.
    Thus the surfaces $\bar S_r$  are supersingular
    for every $r=2^n$,  with Picard number     $16r^2+4r+2$.
\end{say}      
    
\begin{exmp}\label{ss.k3.exmp}
The surface $S_1$ is  a K3 surface with $A_{15}+5A_1$ singularities. 
Resolving the $A_{15}$ singularity, the birational transforms of the $B_i$ become smooth rational curves, and the dual configuration is
    $$
    \begin{array}{ccccccccccccc}
      &&&& \bar B_1 &&&& \bar B_2 &&&& \\
       &&&& \vert &&&& \vert &&&& \\
      C_1 &-& C_2 &-& C_3 &-& \cdots &-& C_{13} &-& C_{14} &-& C_{15}
    \end{array}
    $$
    From this we compute that $(B_1\cdot B_2)=\frac{9}{16}$ and
    $(B_1^2)=( B_2^2)=\frac{7}{16}$.

    Contracting all but the curve $C_8$ gives a supersingular K3 surface $S'_1$ with singularities $2E_8+5A_1$, and Picard number  $1$. The image of $C_8$ is a rational curve $C'_8$ with 2 cusps at the $E_8$-points,  and  $(C'_8\cdot C'_8)=2$.
    Thus $S'_1$ is a  K3 surface of degree 2.
    This is the 4th surface on the list  of \cite{MR2059747}
    in characteristic 2.
\end{exmp}

      \begin{ack}
    I thank M.~Sch\"utt and I.~Shimada  for numerous e-mails explaining supersingular K3 surfaces, and S.~Orevkov and M.~Zaidenberg for references.
      Partial  financial support    was provided  by the Simons Foundation   (grant number SFI-MPS-MOV-00006719-02).
    \end{ack}


\def\cprime{$'$} \def\cprime{$'$} \def\cprime{$'$} \def\cprime{$'$}
  \def\cprime{$'$} \def\dbar{\leavevmode\hbox to 0pt{\hskip.2ex
  \accent"16\hss}d} \def\cprime{$'$} \def\cprime{$'$}
  \def\polhk#1{\setbox0=\hbox{#1}{\ooalign{\hidewidth
  \lower1.5ex\hbox{`}\hidewidth\crcr\unhbox0}}} \def\cprime{$'$}
  \def\cprime{$'$} \def\cprime{$'$} \def\cprime{$'$}
  \def\polhk#1{\setbox0=\hbox{#1}{\ooalign{\hidewidth
  \lower1.5ex\hbox{`}\hidewidth\crcr\unhbox0}}} \def\cdprime{$''$}
  \def\cprime{$'$} \def\cprime{$'$} \def\cprime{$'$} \def\cprime{$'$}
\providecommand{\bysame}{\leavevmode\hbox to3em{\hrulefill}\thinspace}
\providecommand{\MR}{\relax\ifhmode\unskip\space\fi MR }
\providecommand{\MRhref}[2]{%
  \href{http://www.ams.org/mathscinet-getitem?mr=#1}{#2}
}
\providecommand{\href}[2]{#2}

\bigskip

  Princeton University, Princeton NJ 08544-1000, 
\email{kollar@math.princeton.edu}

    \end{document}